\documentclass{amsart}

\usepackage{amssymb,amscd,amsthm,latexsym, amsmath}
\usepackage[all]{xy}

\theoremstyle{plain}
\newtheorem{Lemma}{Lemma}[section]
\newtheorem{theorem}[Lemma]{Theorem}
\newtheorem{lemma}[Lemma]{Lemma}
\newtheorem{proposition}[Lemma]{Proposition}

\theoremstyle{definition}
\newtheorem{definition}[Lemma]{Definition}

\newtheorem{examples}[Lemma]{Examples}

\theoremstyle{remark}

\def\into{ \rightarrowtail }
\def\onto{ \twoheadrightarrow }
\def\splito{ \rightleftarrows }
\newdir{ >}{{}*!/-8pt/@{>}}

\def\CC{ \mathbb{C} }

\def\cleft{\hbox{[\kern-.16em\hbox{[}}}
\def\cright{\hbox{]\kern-.16em\hbox{]}}}

\newcommand{\Ker}{ \ensuremath{\mathrm{Ker}} }

\newcommand{\Rng}{\mathsf{Rng}}

\newcommand{\Gp}{\mathsf{Gp}}

\newcommand{\DiGp}{\mathsf{DiGp}}
\newcommand{\Id}{\mathsf{Id}}

\DeclareMathOperator{\Aut}{Aut}

\newcommand{\SkB}{\mathsf{SkB}}

\numberwithin{equation}{section}

\begin{document}
	
\title{Split epimorphisms and Baer sums\\ of left skew braces}

\author{Dominique Bourn}
\address{Univ. Littoral C\^ote d'Opale, UR 2597, LMPA, Laboratoire de Math\'ematiques Pures et Appliqu\'ees Joseph Liouville, F-62100 Calais, France}
\email{bourn@univ-littoral.fr}


\begin{abstract} We investigate the split epimorphisms in the categories of digroups and left skew braces. We show that, unlike the category $\DiGp$ of digroups, the category $\SkB$ of left skew braces is strongly protomodular. From that, we describe the expected Baer sums of exact sequences of left skew braces  with abelian kernel.\end{abstract}

\maketitle

\section*{Introduction: left skew braces and digroups}

Braces are algebraic structures which were introduced by Rump \cite{Rump} (2007) as producing set-theoretical solutions of the Yang-Baxter equation in response to a general incitement of Drinfeld \cite{Drinfeld} to investigate this equation from a set-theoretical perspective. Later on, Guarnieri and Vendramin generalized this notion with the structure of left skew brace \cite{GV} (2017) which again generates solutions of the Yang-Baxter equations:\\
\noindent\textbf{Definition} 
\emph{A left skew brace is a set $X$ endowed with two group structures $(X,*,\circ)$ subject to a unique axiom: $a\circ (b*c)=(a\circ b)*a^{-*}*(a\circ c)\; (1)$, where $a^{-*}$ denotes the inverse for the law $*$.}

The simplest examples are the following ones: starting with any group $(G,*)$,\\ then $(G,*,*)$ and $(G,*,*^{op})$ are skew braces.

The $*$-multiplication by $a^{-*}$ on the left hand side immediately shows that this axiom is equivalent to the following one: the mapping $\lambda^X_a: X\to X;\; b\mapsto a^{-*}*(a\circ b)$\linebreak is a $*$-homomorphism for any $a\in X$. A major result is that $(1)$ is equivalent to $(2): \lambda^X_{a\circ b}= 
\lambda^X_a\lambda^X_b$ for any $(a,b)\in X\times X$. Furthermore, in a skew brace, the two units necessarily coincide. Taking $a=1^{\circ}=c$ in (1), we get: $b*1^{\circ}=b$, whence $1^{\circ}=1^*$. 

Accordingly, any skew brace is a \emph{digroup}, namely a set $X$ endowed with two group structures $(X,*,\circ)$ only coinciding on the unit, a structure which emerged during discussions between the author and G. Janelidze, and was introduced in \cite{strong} (2000). So, the category $\SkB$ of left skew braces immediately appears as a fully faithful subcategory $\SkB\hookrightarrow \DiGp$ of the category $\DiGp$ of digroups which is clearly stable under products, subobjects and quotient:
namely a Birkhoff subcategory of $\DiGp$. In this way, the category $\SkB$ inherits the structure of \emph{pointed protomodular} category \cite{B0} which is fulfilled by $\DiGp$ and which immediately makes $\SkB$ sharing with the category $\Gp$ of groups all the classical homological lemmas, see \cite{BB}. A first step of investigation of $\SkB$ in this direction was made in \cite{BFP}.

To define protomodularity, only one categorical tool is needed: \emph{the covering pairs of subobjects}, namely pairs $(u,v)$ of subobjects such that any subobject $w: W\into X$ containing them:
$$
\xymatrix@=4pt{  
	&&&	W \ar@{ >->}[dd]^w\\
	&&&&&\\
	U \ar@{ >->}[rrr]_u \ar@{.>}[rrruu]&&& X &&& \ar@{ >->}[lll]^v \ar@{ .>}[llluu] V	  
}
$$
is an isomorphism; i.e. pairs  $(u,v)$ of subobjects such that $U\vee V=X$.\\
\noindent\textbf{Definition} \cite{B0} 
\emph{A pointed category $\CC$ is protomodular when any split epimorphism
	$(f,s)$ is such that the pair $(\ker f,s)$ of subobjects is a covering pair:
	$$
	\xymatrix@=4pt{  
		\Ker f	\ar@{ >->}[rrr]^{\ker f} \ar[dd]&&&	X \ar[dd]_f\\
		&&&&&\\
		1 \ar@{ >->}[rrr]_{0_Y}&&& Y \ar@<-1ex>@{ >->}[uu]_s	  
	}
	$$}\\
Pointed protomodular varieties of universal algebras are characterized in \cite{BJ}.
Pointed protomodularity is a categorical   context in which many basic properties of the paradigmatic example of the category $\Gp$ of groups intrinsiquely hold:\\
- 1) a morphism $f: X\to Y$ is a monomorphism if and only if $\Ker f=1$;\\	
- 2) a regular epimorphism $f: X\onto Y$ is the cokernel of its kernel;\\	
- 3) there is an intrinsic notion of normal subobject; whence an intrinsic notion of exact sequence;\\
-4) there is an intrinsic notion of commutation of pairs $(u: U\into X,v: V\into X)$ of subobjects (denoted by $[u,v]=0$) when there is a (unique possible) homomorphism $+$ in $\CC$ making the following diagram commute:
$$
\xymatrix@=4pt{  
	&&&	U\times V \ar[dd]^+\\
	&&&&&\\
	U \ar@{ >->}[rrr]_u \ar@{ >->}[rrruu]^{\iota_U} &&& X &&& \ar@{ >->}[lll]^v \ar@{ >->}[llluu]_{\iota_V} V	  
}
$$
and consequenly an intrinsic notion of abelian objects $X$ when $[1_X,1_X]=0$;\\
-5) any reflexive relation is an equivalence relation (i.e. protomodularity implies Mal'tsevness).

The abelian objects in $\DiGp$ and $\SkB$ are of the type $(A,+,+)$ where $+$ is commutative. In \cite{strong}, the normal monomorphisms in $\DiGp$ where characterized as 
:\\
\emph{those subobjects $(K,*,\circ)\into (G,*,\circ)$ in $\DiGp$ which are such that:\\
	i) $(K,*)\into (G,*)$ is a normal subgroup;\\
	ii) $(K,\circ)\into (G,\circ)$ is a normal subgroup;\\
	iii) $y*x^{-*}\in K \iff y\circ x^{-\circ}, \;	\forall (x,y) \in G\times G$.
}. 

On their side, Garnieri and Vendramin identified in $\SkB$ a class of subobjects  which they called \emph{ideals}, which actually coincides with the class of normal subobjects in $\SkB$, namely:\\
\emph{those subobjects $(K,*,\circ)\into (G,*,\circ)$ in $\SkB$ which are such that:\\
i)  $(K,\circ)\into (G,\circ)$ is a normal subgroup;\\
ii) $K*a=a*K,\;\forall a\in G$;\\
iii) $\lambda_a(K)\subset K,\;\forall a\in G$.
}\\
Since $\SkB$ is a Birkhoff subcategory of $\DiGp$, any normal subobject  in $\DiGp$ of a left skew brace $(G,*,\circ)$ is an ideal. 
	
Of course, there is a strong relationship between protomodular and additive categories:\\
\emph{
	For any finitely complete pointed category $\CC$, the following conditions are equivalent:\\
	1) $\CC$ is additive;\\
	2) $\CC$ is protomodular and any object is abelian;\\
	3) $\CC$ is protomodular and any subobject is normal.}

In this article, we shall begin by focusing our attention on the characterization of the split epimorphisms in the categories $\DiGp$ and $\SkB$ since they are structural tools in a protomodular context. In spite of the apparent easiness of the question, it hides some unexpected difficulties. 

A split epimorphism $(f,s): (X,*,\circ) \splito (Y,*,\circ)$ in $\DiGp$ is given by a pair of split epimorphisms $(f,s): (X,*) \splito (Y,*), \; (f,s): (X,\circ) \splito (Y,\circ)$ in $\Gp$. And consequently by a pair of group homomorphims (= group actions) $\phi_*: (Y,*)\to (\Aut(K,*),.)$, $\phi_{\circ}: (Y,\circ)\to (\Aut(K,\circ),.)$, where $(K,*,\circ)$ denotes the kernel of $f$ in $\DiGp$. In turn, this pair of actions, through the semidirect products of groups produces a split epimorphism $(p_Y,\iota_Y): (Y\times K,\ltimes_{\phi_\star}, \ltimes_{\phi_{\circ}})\splito (Y,*,\circ)$ in $\DiGp$. However these two split epimorphism \emph{are not isomorphic} in general inside $\DiGp$ since the isomorphisms $(X,*)\to (Y\times K,\ltimes_{\phi_\star})$ and $(X,\circ)\to (Y\times K,\ltimes_{\phi_\circ})$ in $\Gp$ are not in general the same. So, our first aim will be to describe how to produce, up to isomorphism, all the split epimorphisms in $\DiGp$, and then to give several characterizations of those split epimorphisms which lie the subcategory $\SkB$. For that we shall introduce the crucial notion of \emph{index} of a split epimorphims in $\DiGp$. The split epimorphisms in $\DiGp$ and $\SkB$ are investigated by A. Facchini and M. Pompili in \cite{FP} as well, but in a different perspective and with different tools, methods and  achievements; we are grateful to them for several discussions concerning these matters.

From that, we shall show that the category $\SkB$ is not only protomodular but also strongly protomodular \cite{strong} as well, a property shared by the categories $\Gp$ of groups, the category $\Rng$ of rings, and any category of $R$-algebras, which raises $\SkB$ to the level of the richest algebraic structures. On the contrary, the category $\DiGp$ is not strongly protomodular \cite{strong}, it was precisely introduced for this purpose.

It is fascinating to think that any construction in the category $\SkB$ (potentially) produces new solutions to the set-theoretical Yang-Baxter equation. So, we shall finish this article by describing the construction of the Baer sums of exact sequences with abelian kernel in $\SkB$. This construction is expected, as shown in \cite{Baer}, but it will be described here 
in details and under the simplification brought by the "Huq=Smith" condition which is satisfied by the category $\SkB$, see \cite{BFP}.

\section{Index of a split epimorphism}\label{index}

Suppose given any split epimorphism $(f,s): (X,1)\splito (Y,1)$ in the category $Set_*$ of pointed sets, and denote by $(K,1)$ its kernel.
\begin{definition}
	An index for the split epimorphism $(f,s)$ is a bijection $\chi: (X,1)\to (X,1)$ satisfying $f\chi=f$, $\chi s=s$ and such that the restriction $\chi/_K$ of $\chi$ to the kernel $K$ of $f$ is $\Id_K$.
\end{definition}
Clearly the indexes of $(f,s)$ determines a subgroup of the symmetric group $S_X$. 
\begin{proposition}
	Given $(Y,1)$ and $(K,1)$ any pair of pointed sets, the indexes of the canonical split epimorphism $(p_Y,\iota_Y): (Y\times K,1) \splito (Y,1)$ are in one to one correspondance with the mappings $\xi\colon (Y,1)\to (S_K,\Id)$ in $Set_*$
\end{proposition}
\proof
The index $\chi$ associated with the mapping $\xi$ is defined by $\chi(y,k)=(y,\xi_y(k))$, and conversely the mapping associated with the index $\chi$ is defined by $\xi_y(k)=p_K.\chi(y,k)$. We shall keep the same word \emph{index} for those mappings $\xi$.
\endproof

Given any split epimorphism $(f,s): (X,\circ)\splito (Y,\circ)$ in the categgory $\Gp$ of groups, we shall denote by $\phi_{\circ}: (Y,\circ) \to \Aut(K,\circ)$, the classical associated group action defined by $k\mapsto s(y)\circ k\circ s(y)^{-\circ}$, and by $\tau_{(f,s)}^{\circ}:(X,1)\to (K,1)$ the mapping in $Set_*$ defined by $\tau_{(f,s)}^{\circ}(x)=x\circ sf(x)^{-\circ}$. Conversely, starting from a group action $\phi: (Y,\circ) \to \Aut(K,\circ)$, we shall denote the split epimorphism obtained by the semi-direct product in the following (non classical) way: $(p_Y,\iota_Y): (Y\times K, \ltimes_\phi) \splito (Y,\circ)$ in order to avoid any confusion arising from the multiple possible group laws existing on the sets $Y$ and $K$.

\medskip

Now, working in the category $\DiGp$ of digroups, and since the lower square in the following diagram is a pullback in $Set_*$, any split epimorphism $(f,s)\colon (X,*,\circ) \splito (Y,*,\circ)$ in $\DiGp$ produces a unique mapping $\chi\colon  (X,1)\to (X,1)$ in $Set_*$ making the following diagram commute, :
$$
\xymatrix@=10pt{ 
	(X,1) \ar[dddr]_{f} \ar[rrrd]^{\tau_{(f,s)}^{*}} \ar@{.>}[rd]_>>>{\chi} &&&& \\
	&	(X,1) \ar[dd]^{f}  \ar[rr]_{\tau_{(f,s)}^{\circ}}  && (K,1) \ar[dd]^{\tau_K}\\
	&&&   \\
	& (Y,1)	 \ar[rr]_{\tau_Y} && 1  
}
$$
This mapping is a bijection, since the upper quadrangle is a pullback as well.
This natural comparison mapping  $\chi$ is defined by $x\mapsto (x*sf(x)^{-*})\circ sf(x)$, and it is clearly an index for the underlying split epimorphism $(f,s)$ in $Set_*$. We shall call it the \emph{index} of the split epimorphism $(f,s)\colon (X,*,\circ) \splito (Y,*,\circ)$ in $\DiGp$.  When $(X,*,\circ)=(X,*,*)$, the index is trivial; when $(X,*,\circ)=(X,*,*^{op})$, we get $\chi(x)=sf(x)*x*sf(x)^{-*}$.

\medskip

We get $\chi(k*s(y))=k\circ s(y)$  for all $(y,k)\in Y\times K$.
In turn, this index $\chi$ for $(f,s)$ produces an index for the split epimorphism $(p_Y,\iota_Y):Y\times K \splito Y$ defined by $\xi^{(f,s)}_y(k)=\chi(k*s(y))*s(y)^{-*}=(k\circ s(y))*s(y)^{-*}$.

\begin{definition}\label{skewi}
	Given any split epimorphism $(f,s)\colon (X,*,\circ) \splito (Y,*,\circ)$ in $\DiGp$, the mapping $\xi^{(f,s)}:(Y,1)\to  (S_K,\Id_K)$ in $Set_*$ is called  the \emph{skewing index} of $(f,s)$.
\end{definition}

\begin{proposition}\label{prop3}
	For any split epimorphism $(f,s)\colon (X,*,\circ) \splito (Y,*,\circ)$ in $\DiGp$  the following conditions are equivalent:\\
	1) $\chi=\Id_X$, i.e. the index is trivial;\\
	2) $x\circ sf(x)^{-\circ}=x*sf(x)^{-*}$ for all $x\in X$;\\
	3) $k\circ s(y)=k*s(y)$ for all $(y,k)\in Y\times K$;\\
	4) the skewing index $\xi^{(f,s)}$ is trivial, namely such that $\xi^{(f,s)}_y=\Id_K$ for any $y\in Y$;\\
	5) we have $\lambda^X_ks=s$ for any $k\in K$.
\end{proposition}
\proof
We get  $\chi=\Id_X$ if and only if $(x*sf(x)^{-*})\circ sf(x)=x$
namely if and only if $x\circ sf(x)^{-\circ}=x*sf(x)^{-*}$.

\noindent We already noticed that:
$\chi(k*s(y))=k\circ s(y)$.
So, $\chi=\Id_X$ holds if and only if $k\circ s(y)=k*s(y)$ for all $(y,k)\in Y\times K$. The points 4) and 5) are then straightforward.
\endproof

\begin{lemma}
	Any split epimorphism $(Y\times K,\ltimes_{\psi_{*}},\ltimes_{\psi_{\circ}})\splito  (Y,*,\circ)$ in $\DiGp$ built from a pair $(\psi_{*}: (Y,*)\to Aut(K,*),\psi_{\circ}: (Y,\circ)\to Aut(K,\circ))$ of group actions  has a trivial index.
\end{lemma}
\proof
We get $(1,k)\ltimes_{\psi_{*}}(y,1)=(y,k)=(1,k)\ltimes_{\psi_{\circ}}(y,1)$, namely 3). 
So, the index $\chi$ is $\Id_{Y\times K}$
\endproof

Actually, this is a characterisation, in the category $\DiGp$, of the split epimorphisms arising from a pair $(\psi_{*},\psi_{\circ})$  of group actions:

\begin{theorem}\label{Id}
	Let $(f,s)\colon  (X,\star,\circ)\splito (Y,\star,\circ)$ be any split epimorphism in $\DiGp$ with kernel $(K,\star,\circ)$, and $(\psi_{\star}, \psi_{\circ})$ be any pair of group actions. Then the two following conditions are equivalent:\\
	1) there exists an isomorphism $\alpha$ in $\DiGp$ making the following diagram commute:
	$$
	\xymatrix@=4pt{ 
		&	(X,\star,\circ) \ar[dd]_{f}  \ar[rr]^<<<<<{\alpha}  && (Y\times K,\ltimes_{\psi_{\star}},\ltimes_{\psi_{\circ}}) \ar[dd]_{p_Y}\\
		&&&   \\
		& (Y,\star,\circ)	\ar@<-1ex>[uu]_{s} \ar[rr]_{Id_Y} && (Y,\star,\circ) \ar@<-1ex>[uu]_{\iota_Y} 
	}
	$$
	2) the index $\chi$ of $(f,s)$ is trivial. 
\end{theorem}
\proof
Suppose $(f,s)\colon  (X,\star,\circ)\splito (Y,\star,\circ)$ is a split epimorphism with $\chi=\Id_X$, namely such that $x\star sf(x)^{-\star}=x\circ sf(x)^{-\circ}$. Then define $\alpha\colon  X\to Y\times K$ by $x \mapsto (f(x),x\star sf(x)^{-\star})=(f(x),x\circ sf(x)^{-\circ})$. By the first member of the equality, we get a group isomorphism $\alpha\colon  (X,\star)\to (Y\times K,\ltimes_{\psi_*})$; and by the second one  a group isomorphism $\alpha\colon  (X,\circ)\to (Y\times K,\ltimes_{\psi_\circ})$.

Conversely, start with a pair $(\alpha,\Id_Y)$ of digroup isomorphisms as above. From $p_Y\alpha=f$, the map $\alpha$ is produced by a pair of mappings $(f,\gamma)$ with a mapping $\gamma\colon  X \to K$ such that its restriction $\gamma_{/K} \colon K\to K$ produces a digroup automorphism  $(K,\star,\circ) \to (K,\star,\circ)$.
From $\alpha s=\iota_Y$, we get $\gamma s(y)=1$ for all $y\in Y$.

\medskip

Now, since $(f,\gamma)\colon  (X,\star)\to (Y\times K,\ltimes_{\psi_{\star}})$ is a group homomorphism, we get for all $(y,k)\in Y\times K$:
$$(f,\gamma)(k\star s(y))=(f,\gamma)(k)\ltimes_{\psi_{\star}} (f,\gamma)s(y)=(1,\gamma(k))\ltimes_{\psi_{\star}}(y,1)
=(y,\gamma(k)).$$
Whence: $(f,\gamma)(x)=(f,\gamma)((x \star sf(x)^{-\star})\star sf(x))
=(f(x),\gamma(x\star sf(x)^{-\star}))$. 

Since $(f,\gamma)\colon  (X,\circ) \to,(Y\times K,\ltimes_{\psi_{\circ}})$ is a group homomorphism as well, we get:
$(f,\gamma)(x)=(f(x),\gamma(x\circ sf(x)^{-\circ}))$. So:
$\gamma(x\star sf(x)^{-\star})=\gamma(x\circ sf(x)^{-\circ})$. And since the restriction of $\gamma$ to $K$ is a bijection, we get  $x\star sf(x)^{-\star}=x\circ sf(x)^{-\circ}$. By Proposition \ref{prop3}, this does  happen if and only if $\chi=\Id_X$.
\endproof
Given any split epimorphism $(f,s)\colon  (X,\star,\circ)\splito (Y,\star,\circ)$ in $\DiGp$, this theorem shows how unexpectedly independent is the index $\chi$  from the associated pair $(\phi_*,\phi_{\circ})$ of group actions producing separately the split epimorhisms $(f,s)\colon  (X,\star)\splito (Y,\star)$ and $(f,s)\colon  (X,\circ)\splito (Y,\circ)$ in $\Gp$. We get some more precision with the following:
\begin{proposition}\label{split}
 	Let $(f,s)\colon  (X,\star,\circ)\splito (Y,\star,\circ)$ be any split epimorphism in $\DiGp$ with skewing index $\xi^{(f,s)}$. We get: $\xi^{(f,s)}_y\phi_{\circ y}=\phi_{* y}\lambda^X_{s(y)}$ for all $y\in Y$.
\end{proposition}
\proof
The two terms are equal to: $(s(y)\circ k)\star s(y)^{-\star}$.
\endproof

\section{Split epimorphisms of digroups}

It remains to produce, up to isomorphism in $\DiGp$, all the split epimorphisms in $\DiGp$ whose index is not trivial.

Let $(f,s)\colon  (X,\star,\circ)\splito (Y,\star,\circ)$ be any split epimorphism in $\DiGp$ with kernel $(K,\star,\circ)$. We get the two canonical  group actions $\phi_\star: (Y,\star)\to \Aut(K,\star)$ and $\phi_\circ: (Y,\circ)\to \Aut(K,\circ)$, and the skewing index $\xi^{(f,s)}:(Y,1)\to (S_K,\Id_K)$ in $Set_*$ as well.

Conversely, let us start with two group actions $\psi_\star: (Y,\star)\to \Aut(K,\star)$, $\psi_\circ: (Y,\circ)\to \Aut(K,\circ)$, and any mapping $\xi:(Y,1)\to (S_K,\Id_K)$ in $Set_*$. Let us denote  by $\Xi: Y\times K \to Y\times K$, the associated permutation on $Y\times K$ given by $(y,k)\mapsto (y,\xi_y(k))$, and by $\circ_\xi$ the transfer of the group law $\ltimes_{\psi_\circ}$ on $Y\times K$ by the permutation  $\Xi$. We get:
$$(y,k)\circ_\xi(y',k')=\bigl(y\circ y', \xi_{y\circ y'}(\xi_y^{-1}(k)\circ \psi_{\circ y}(\xi_{y'}^{-1}(k')))\bigr) $$ 
\begin{proposition}\label{unit}
	The split epimorphism $(p_Y,\iota_Y): (Y\times K,\ltimes_{\psi_\star},\circ _\xi)\splito (Y,\star,\circ)$ is a split epimorphism in  $\DiGp$ whose kernel is $(K,\star,\circ)$, whose skewing index is $\xi$, and such that $\phi_{\ltimes_{\psi_\star}}=\psi_\star$, $\phi_{\circ_\xi}=\psi_\circ$.
\end{proposition}
\proof
The proof is straightforward. We get the following table:
\begin{center}
	\begin{tabular}{ c | c }
		
		$(1,k)\circ_\xi(1,k')=(1,k\circ k')$ & $(y,1)\circ_\xi(y',1)=(y\circ y',1)$  \\
		$ (1,k)\circ_\xi(y,1)=(y,\xi_y(k))$ &  $\;\;\;\;\;(y,1)\circ_\xi(1,k)= (y,\xi_y\psi_{\circ y}(k))$ \\
		
	\end{tabular}
\end{center}
So, the split epimorphism in question lies in $\DiGp$, and its kernel is $(K,\star,\circ)$. We classically get $\phi_{\ltimes_{\psi_\star}}=\psi_\star$. Now:\\ 
$(y,1)\circ_\xi(1,k)\circ_\xi(y^{-\circ},1)=(y,1)\circ_\xi(y^{-\circ},\xi_{y^{-\circ}}(k))=(1,\psi_{\circ y}\xi_{y^{-\circ}}^{-1}\xi_{y^{-\circ}}(k))=(1,\psi_{\circ y}(k))$
which shows that $\phi_{\circ_\xi}=\psi_\circ$. Finally the index $\chi$ of $(p_Y,\iota_Y)$ is defined by:\\ $\chi(y,k)=((y,k)\ltimes_{\psi_\star}(y^{-\star},1))\circ_\xi(y,1)=(1,k)\circ_\xi(y,1)=(y,\xi_y(k))$, which shows that the skewing index of $(p_Y,\iota_Y)$ is $\xi$. 
\endproof

Conversely we get:
\begin{proposition}
	Let  $(f,s)\colon  (X,\star,\circ)\splito (Y,\star,\circ)$ be any split epimorphism in $\DiGp$ with kernel $(K,\star,\circ)$ and skewing index $\xi$.  Then $(f,s)$ is isomorphic to $(p_Y,\iota_Y): (Y\times K,\ltimes_{\phi_\star},\circ _\xi)\splito (Y,\star,\circ)$ in $\DiGp$.
\end{proposition}
\proof
Classically, the mapping $\theta$ defined by $(y,k)\mapsto k\star s(y)$ determines a group isomorphism $\theta: (Y\times K,\ltimes_{\phi_\star}) \to (X,\star)$.

The skewing index $\xi_y$ of $(f,s)$, see Definition \ref{skewi}, is defined by $k\mapsto (k\circ s(y))*s(y)^{-*}$. So, $(\xi_y)^{-1}(y)=(k*s(y))\circ s(y)^{-\circ}$. From that, we get:\\
$\xi_y^{-1}(k)\circ \phi_{\circ y}(\xi_{y'}^{-1}(k'))=(k*s(y))\circ s(y)^{-\circ} \circ s(y) \circ (k'*s(y'))\circ s(y')^{-\circ}\circ s(y)^{-\circ}$\\
$= (k*s(y))\circ (k'*s(y'))\circ s(y\circ y')^{-\circ}$. So: $$(y,k)\circ_\xi(y',k')=\bigl(y\circ y',[(k\star s(y))\circ (k'\star s(y'))]*s(y\circ y')^{-*}\bigr)$$ 
Then $\theta((y,k)\circ_\xi(y',k'))=(k\star s(y))\circ (k'\star s(y'))=\theta(y,k)\circ \theta(y',k')$. Accordingly, $\theta:(Y\times K,\circ_\xi)\to(X,\circ)$ is a group homomorphism as well.
\endproof

\section{Split epimorphisms of skew braces}

We have now to identify the split epimorphisms  $(f,s)\colon  (X,\star,\circ)\splito (Y,\star,\circ)$ of skew braces among the split epimorphisms of digroups. Of course a necessary and sufficient condition is that $(X,\star,\circ)$ is a skew brace, namely that $\lambda^X_x:X\to X$ is a $*$-homomorphism for any $x\in X$. However, beyond other characterizations, our investigation will produce easy obstruction conditions and sufficient conditions which will allow us to build up some examples. 

\subsection{The general case}

A first simplification step is the following:

\begin{lemma}\label{prop4}
	Let $(f,s): (X,*,\circ) \splito (Y,*,\circ)$ be any split epimorphism in $\DiGp$. The following conditions are equivalent:\\
	$\alpha$) $\lambda^X_{s(y)}:(X,*)\to (X,*)$ is a group homomorphism for any $y\in Y$;\\
	$\beta$) i) $(Y,*,\circ)$ is a skew brace;\\
   ii) the restriction of $\lambda^X_{s(y)}$ to $K$ is a group automorphism $(K,*)\to (K,*)$ for all $y\in Y$;\\
   iii) $\lambda^X_{s(\gamma)}\phi_{* y}(k)*s\lambda^Y_\gamma(y)=s\lambda^Y_\gamma(y)*\lambda^X_{s(\gamma)}(k)$ for all $(\gamma,y,k)\in Y\times Y \times K$.
	
\noindent	When any the the previous conditions hold, we get $\lambda^X_{s(y)\circ k}=\lambda^X_{s(y)}\lambda^X_{k}$ for all $(y,k)\in Y\times K$, and $\lambda^X_{s(y\circ y')}=\lambda^X_{s(y)}\lambda^X_{s(y')}$ for all $(y,y')\in Y\times Y$. 
\end{lemma}
\proof
1) We need the "restriction" of $\lambda^X_{s(y)}$ to $Y$ which is $\lambda^X_{s(y)}s=s.\lambda^Y_{y}$, and the restriction of $\lambda^X_{s(y)}$ to $K$. When $\lambda^X_{s(y)}$ is a group homomorphism, both are group homomorphisms.

Conversely saying that  the "restriction" $s.\lambda^Y_{y}$ is a group homomorphism, is equivalent ($s$ being injective) to saying that $\lambda^Y_{y}$ is a group homomorphism $(Y,*)\to (Y,*)$ for all $y\in Y$; namely that $(Y,*,\circ)$ is a skew brace.
So, under conditions $\beta$i and $\beta$ii we can consider the following diagram in $\Gp$:
	$$\xymatrix@=7pt{ 
	&&& (X,\star)\\
	(K,\star) \ar[dd]  \ar@{ >->}[rr]_{k_f} \ar@{ >->}[rrru]^{\lambda^X_{s(y)}}  && (X,\star) \ar@{.>}[ur] \ar[dd]_{f}\\
	&&&   \\
	1	\ar@<-1ex>[uu] \ar@{ >->}[rr]
	&& (Y,\star) \ar@<-1ex>[uu]_{s} \ar@<-1ex>[uuur]_{s\lambda^Y_y}
}$$
Applying Lemma \ref{joint} (below) determines condition $\beta$iii.

\smallskip

\noindent 2) When $\lambda^X_{s(y)}: (X,*) \to (X,*)$ is a group homomorphism, then $\lambda^X_{s(y)}(k^{-*})=(\lambda^X_{s(y)}(k))^{-*}$. Whence: $(s(y)\circ k)^{-*}=s(y)^{-*}*(s(y)\circ k^{-*})*s(y)^{-*}$, and:\\
$\lambda^X_{s(y)\circ k}(x)=(s(y)\circ k)^{-*}*(s(y)\circ k \circ x)=s(y)^{-*}*(s(y)\circ k^{-*})*s(y)^{-*}*(s(y)\circ k \circ x)$
$=\lambda^X_{s(y)}(k^{-*})*\lambda^X_{s(y)}(k\circ x)=\lambda^X_{s(y)}(k^{-*}*(k\circ x))=\lambda^X_{s(y)}\lambda^X_{k}(x)$.

\smallskip
When $(Y,*,\circ)$ is a skew brace, then $\lambda^Y_{y}(y'^{-*})=\bigl(\lambda^Y_{y}(y')\bigr)^{-*}$.\\  Whence: $(s(y)\circ s(y'))^{-*}=s(y)^{-*}*(s(y)\circ s(y')^{-*})*s(y)^{-*}$, and:\\
$\lambda^X_{s(y)\circ s(y')}(x)=(s(y)\circ s(y'))^{-*}*(s(y)\circ s(y') \circ x)$\\$=s(y)^{-*}*(s(y)\circ s(y')^{-*})*s(y)^{-*}*(s(y)\circ s(y') \circ x)$
$=\lambda^X_{s(y)}(s(y')^{-*})*\lambda^X_{s(y)}(s(y')\circ x)$\\$=\lambda^X_{s(y)}(s(y')^{-*}*(s(y')\circ x))=\lambda^X_{s(y)}\lambda^X_{s(y')}(x)$. 
\endproof
\begin{theorem}
	Let $(f,s): (X,*,\circ) \splito (Y,*,\circ)$ be any split epimorphism in $\DiGp$. The following conditions are equivalent:\\
	1) $(f,s): (X,*,\circ) \splito (Y,*,\circ)$ is a split epimorphism in $\SkB$;\\
	2) (i) $\lambda^X_{s(y)}:(X,*)\to (X,*)$ and (ii) $\lambda^X_k: (X,*)\to (X,*)$ are group homomorphisms for any $(y,k)\in Y\times K$.
\end{theorem}
\proof
The split epimorphism $(f,s)$ is in $\SkB$ if and only if $\lambda^X_x$ is a $*$-homomorphism for any $x\in X$. So that 1) implies 2) is straightforward. Suppose 2). Then, according to the previous proposition, 2(i) implies that  $\lambda^X_{s(y)\circ k}=\lambda^X_{s(y)}\lambda^X_{k}$ for all $(y,k)\in Y\times K$; with 2(ii) we can check that that any  $\lambda^X_x=\lambda^X_{s(y)\circ k}: (X,*)\to (X,*)$ is a group homomorphism. Accordingly $(X,*,\circ)$ is a left skew brace. 
\endproof
The previous lemma produces a symmetric assertion:
\begin{lemma}\label{prop8}
	Let $(f,s): (X,*,\circ) \splito (Y,*,\circ)$ be any split epimorphism in $\DiGp$. The following conditions are equivalent:\\
	$\alpha$) $\lambda^X_{k}:(X,*)\to (X,*)$ is a group homomorphism for any $k\in K$;\\
	$\beta$) i) $(K,*,\circ)$ is a skew brace;\\
	ii) the "restriction" $\lambda^X_{k}s$ of $\lambda^X_{k}$ to $Y$ is a homomorphic section $(Y,*)\to (X,*)$ of $f$ for all $k\in K$;\\
	iii) $\lambda^X_{k}\phi_{* y}(\bar k)*\lambda^X_k(s(y))=\lambda^X_k(s(y))*\lambda^X_k(\bar k)$ for all $(y,k,\bar k)\in Y\times K \times K$.
	
	\noindent	When any the the previous conditions hold, we get $\lambda^X_{k\circ s(y)}=\lambda^X_{k}\lambda^X_{s(y)}$ for all $(y,k)\in Y\times K$, and $\lambda^X_{k\circ k'}=\lambda^X_{k}\lambda^X_{k'}$ for all $(k,k')\in K\times K$. 
\end{lemma}
\proof
1) The restriction of $\lambda^X_k$ to $K$ is $\lambda^X_kk_f=k_f\lambda^K_k$. Saying it is a group homomorphism is saying ($k_f$ being injective) that $\lambda^K_k$ is a group homomorphism $(K,*)\to (K,*)$ for all $k\in K$, namely that $(K,*,\circ)$ is a skew brace. Then,
this time, we apply Lemma \ref{joint} to the following diagram in $\Gp$:
	$$\xymatrix@=7pt{ 
	&&& (X,\star)\\
	(K,\star) \ar[dd]  \ar@{ >->}[rr]_{k_f} \ar@{ >->}[rrru]^{k_f\lambda^Y_k}  && (X,\star) \ar@{.>}[ur] \ar[dd]_{f}\\
	&&&   \\
	1	\ar@<-1ex>[uu] \ar@{ >->}[rr]
	&& (Y,\star) \ar@<-1ex>[uu]_{s} \ar@<-1ex>[uuur]_{\lambda^X_ks}
}$$

\noindent 2) This times we use the two following facts: $\lambda^X_{k}(s(y)^{-*})=(\lambda^X_{k}(s(y)))^{-*}$ and $\lambda^X_{k}(k'^{-*})=(\lambda^X_{k}(k'))^{-*}$.
\endproof

We are now in position to assert:

\begin{proposition}\label{final}
	The split epimorphism $(f,s): (X,*,\circ) \splito (Y,*,\circ)$ with kernel $(K,*,\circ)$ in $\DiGp$ belongs to $\SkB$ if and only if:\\
	i) $(Y,*,\circ)$ and $(K,*,\circ)$ are skew braces;\\
	ii) - the restriction of $\lambda^X_{s(y)}$ to $K$ is a group automorphism $(K,*)\to (K,*)$ for all $y\in Y$, and\\
	- the "restriction" $\lambda^X_{k}s$ of $\lambda^X_{k}$ to $Y$ is a homomorphic section $(Y,*)\to (X,*)$ of $f$ for all $k\in K$;\\
	iii) $\lambda^K_{k}\phi_{* y}(\bar k)*\lambda^X_ks(y)=\lambda^X_ks(y)*\lambda^K_k(\bar k)$, for all $(y,k,\bar k)\in Y\times K\times K$
	
	$\lambda^X_{s(\gamma)}\phi_{* y}(k)*s\lambda^Y_\gamma(y)=s\lambda^Y_\gamma(y)*\lambda^X_{s(\gamma)}(k)$, for all $(\gamma,y,k)\in Y\times Y \times K$.
\end{proposition}

\subsection{Split epimorphisms with trivial index in $\SkB$}

Let us start now with  any pair $\psi_{*}: (Y,*)\to \Aut(K,*)$ and $\psi_{\circ }: (Y,\circ)\to \Aut(K,\circ)$ of group actions. When does the split epimorphism $(p_Y,\iota_Y): (Y\times K,\ltimes_{\psi_\star},\ltimes_{\psi_\circ}) \splito (Y,*,\circ )$ in $\DiGp$ belong to $\SkB$?
\begin{theorem}
	Given any pair $(Y,*,\circ),(K,*,\circ)$ of skew braces, and any pair $\psi_{*}: (Y,*)\to \Aut(K,*)$ and $\psi_{\circ }: (Y,\circ)\to \Aut(K,\circ)$ of group actions, the split epimorphism $(p_Y,\iota_Y): (Y\times K,\ltimes_{\psi_\star},\ltimes_{\psi_\circ}) \splito (Y,*,\circ )$ belongs to $\SkB$ if and only if:\\
	1) $\psi_{\circ y}$ is a $*$-automorphism for all $y\in Y$;\\
	2) $\psi_{* y}\lambda^K_k=\lambda^K_k\psi_{* y}$ for all $(y,k)\in Y\times K$;\\
	3) $\psi_{* (\gamma\circ y)*\gamma^{-*}}=\psi_{\circ \gamma}\psi_{* y}\psi_{\circ \gamma}^{-1}$ for all $(\gamma,y)\in Y\times Y$.
\end{theorem}
\proof
Let us make explicit the mapping $\lambda^{Y\times K}_{(y,k)}$. We get:\\
$(y,k)^{-\ltimes_{\psi_\star}}\ltimes_{\psi_\star}\bigl((y,k)\ltimes_{\psi_{\circ}}(\bar y, \bar k)\bigr)=(y^{-*},\psi_{* y^{-*}}(k^{-*})){\ltimes_{\psi_\star}}(y\circ \bar y,k\circ \psi_{\circ y}(\bar k)))$\\
$=(y^{-*}*(y\circ \bar y),\psi_{* y^{-*}}(k^{-*})*\psi_{* y^{-*}}(k\circ \psi_{\circ y}(\bar k))=(\lambda^Y_y(\bar y),\psi_{* y^{-*}}(k^{-*}*(k\circ \psi_{\circ y}(\bar k))))$, namely, $\lambda^{Y\times K}_{(y,k)}(\bar y,\bar k)=(\lambda^Y_y(\bar y),\psi_{* y^{-*}}(\lambda^K_k(\psi_{\circ y}(\bar k))))$. In particular, we get:\\
$\lambda_{(y,1)}(\bar y,\bar  k)=(\lambda^Y_y(\bar y),\psi_{* y^{-*}}\psi_{\circ y}(\bar k))$ and $\lambda_{(1,k)}(\bar y,\bar  k)=(\bar y,\lambda^K_k(\bar k))$.

So, the restriction of $\lambda_{(y,1)}$ to $K$ is $k\mapsto \psi_{* y^{-*}}\psi_{\circ y}(k)$, while
the "restriction" of $\lambda_{(1,k)}$ to $Y$ is $y\mapsto (y,1)=\iota_Y(y)$ for any $k\in K$.\\
And, as well:
$\lambda_{(y,1)}(\bar y,1)=(\lambda^Y_y(\bar y),1)$ while $\lambda_{(1,k)}(1,\bar  k)=(1,\lambda^K_k(\bar k))$. 

Now, according to the previous proposition, the split epimorphism in question belongs to $\SkB$ if and only if:\\
ii) - $k\mapsto \psi_{* y^{-*}}\psi_{\circ y}(k)$ is a group automorphism $(K,*)\to (K,*)$ for all $y\in Y$, which is equivalent to $\psi_{\circ y}(k)$ is a $*$-automorphism, since so is $\psi_{* y^{-*}}$, which, in turn, is equivalent  to condition 1) in our assertion;\\
- $y\mapsto (y,1)=\iota_Y(y)$ is a homomorphic section $(Y,*)\to (X,*)$ of $p_Y$, which is trivial;\\
iii) $(1,\lambda^K_{k}\phi_{* y}(\bar k))\ltimes_{\psi_*}(y,1)=(y,1)\ltimes_{\psi_*}(1,\lambda^K_k(\bar k))$, for all $(y,k,\bar k)\in Y\times K\times K$, namely: $(y,\lambda^K\psi_{* y}(\bar k))=(y,\psi_{* y}\lambda^K_k(\bar k))$, i.e condition 2) in our assertion; and\\
$(1,\psi_{* \gamma^{-*}}\psi_{\circ \gamma}\psi_{* y}(k))\ltimes_{\psi_*}(\lambda^Y_\gamma (y),1)=(\lambda^Y_\gamma (y),1)\ltimes_{\psi_*}(1,\psi_{* \gamma^{-*}}\psi_{\circ \gamma}(k))$, namely:\\
$(\lambda^Y_\gamma (y),\psi_{* \gamma^{-*}}\psi_{\circ \gamma}\psi_{* y}(k))=(\lambda^Y_\gamma (y),\psi_{* \lambda^Y_\gamma(y)}(\psi_{* \gamma^{-*}}\psi_{\circ \gamma}(k))$.\\ This condition is equivalent to:
$\psi_{* \gamma^{-*}}\psi_{\circ \gamma}\psi_{* y}(k)=\psi_{* \gamma^{-*}*(\gamma\circ y)}\psi_{* \gamma^{-*}}\psi_{\circ \gamma}(k)$, i.e:\\
$\psi_{\circ \gamma}\psi_{* y}(k)=\psi_{* (\gamma\circ y)}\psi_{* \gamma^{-*}}\psi_{\circ \gamma}(k)=\psi_{* (\gamma\circ y)*\gamma^{-*}}\psi_{\circ \gamma}(k)$, which, in turn, is equivalent to condition 3).
\endproof
\begin{examples}
	1) Let $(Y,\circ)$ and $(K,\circ)$ be two groups; and let $\psi_\circ$ be the trivial action $\check 1: (Y,\circ)\to \Aut(K,\circ)$; it produces the split epimorphism $(p_Y,\iota_Y)$:\linebreak $(Y\times K, \times_\circ) \splito (Y,\circ)$, where $\times_\circ$ denotes the law of the cartesian product $(Y,\circ)\times(K,\circ)$. Since $\check 1_y=\Id_K$ for any $y\in Y$, the condition 1) of the above theorem trivially holds for any group action $\psi_*: (Y,*)\to \Aut(K,*)$; and the condition 3) becomes $\psi_{* (\gamma\circ y)*\gamma^{-*}}=\psi_{* y}$, namely $\psi_{* (\gamma\circ y)}=\psi_{* y}\psi_{* \gamma}=\psi_{*( y*\gamma)}$. This holds for any $\psi_*$ as soon as $\circ=*^{op}$ on $Y$. Consider the split epimorphism $(p_Y,\iota_Y): (Y\times K,\ltimes_{\psi_{*}},\times_{*^{op}})\splito (Y,*,*^{op})$ in $\DiGp$; it kernel is $(K,*,*^{op})$. Now $\lambda^K_k(\bar k)=k^{-*}*\bar k*k$, and the condition 2) becomes $\psi_{* y}(k^{-*}*\bar k*k)=k^{-*}*\psi_{* y}(\bar k)*k$, which holds as soon as $*$ is commutative on $K$, namely $*=*^{op}$ on $K$. 
	
	Accordingly, given any group action $\psi: (Y,*)\to \Aut(A,+)$, the split epimorphism 
	$$(p_Y,\iota_Y): (Y\times A,\ltimes_{\psi},\times_{*^{op}}) \splito (Y,*,*^{op})$$
	where $\times_{*^{op}}$ denotes the law on the cartesian product $(Y,*^{op
	})\times (A,+)$, belongs to $\SkB$ and has an abelian kernel $(A,+,+)$. 
	
	2) Let $(Y,*,\circ)$ be a skew brace and $\psi: (Y,\circ)\to \Aut(A,+)$ any group action. Then the split epimorphism in $\DiGp$:
	$$(p_Y,\iota_Y): (Y\times A,\times, \times_{\psi}) \splito (Y,*,\circ)$$
	where $\times$ denotes the law on the cartesian product $(Y,*)\times (A,+)$
	belongs to $\SkB$ and has an abelian kernel $(A,+,+)$.
	\proof
	This times, it is the group action $\psi_*$ which is trivial. So, 1) $\psi_{y}$ is a $+$-automorphism for any $y\in Y$. 2) $\lambda^A_a
=\Id_A$	commutes with any automorphism. And condition 3) becomes $\psi_{\gamma}\psi_{\gamma}^{-1}=\Id_A$ which is trivial.\endproof

3) So, starting with two abelian groups $(B,+)$ and $(A,+)$ and a group action $\psi: (B,+)\to \Aut(A,+)$, both split epimorphisms $(p_{B},\iota_{B}): (B\times A, \times,\ltimes_{\psi})\splito (B,+,+)$ and $(p_{B},\iota_{B}): (B\times A,\ltimes_{\psi},\times)\splito (B,+,+)$ lie in $\SkB$ and has abelian kernel $(A,+,+)$.
\end{examples}

\section{Strong protomodularity}

\subsection{Fibration of points and protomodularity}

We recalled the definition of protomodularity in the introduction. The fact that the pair $(k_f,s)$ is jointly strongly epic for any split epimorphism $(f,s)$ produces the  following useful result in the category $\Gp$:

\begin{lemma}\label{joint}
	Let $(X,\star)\splito (Y,\star)$ be any split epimorphism in $\Gp$, and $(K,\star)$ the kernel of $f$:
	$$
	\xymatrix@=7pt{ 
		&&& (D,\star)\\
		(K,\star) \ar[dd]  \ar@{ >->}[rr]_{k_f} \ar@{ >->}[rrru]^{g_K}  && (X,\star) \ar@{.>}[ur]_g \ar[dd]_{f}\\
		&&&   \\
		1	\ar@<-1ex>[uu] \ar@{ >->}[rr]
		&& (Y,\star) \ar@<-1ex>[uu]_{s} \ar@<-1ex>[uuur]_{g_Y}
	}$$
	Then a pair $(g_K,g_Y)$ of group homomorphisms determines a (unique possible) factorization $g: X\to D$ in $\Gp$ if and only if:\\
	$g_K\phi_{\star y}(k)\star g_Y(y)=g_Y(y)\star  g_K(k)$ for all $(y,k)\in Y\times K$,\\ or equivalently: $g_K\phi_{\star y}=\phi_{\star g_Y(y)}g_K$ for all $y\in Y$. 
\end{lemma}
\proof
The commutation of the triangles implies that we have necessarily:\\ $g(k\star s(y))=g_K(k)\star g_Y(y)$. It remains to show that this definition determines a group homomorphism. From $k\star s(y) \star k' \star s(y')=k\star \phi_{\star y}(k') \star s(y\star y')$, this is equivalent to:
$g_K(k)\star g_Y(y) \star g_K(k') \star g_Y(y')=g_K(k\star \phi_{\star y}(k')) \star g_Y(y\star y')$ which, in turn, is equivalent to: $g_Y(y) \star g_K(k') =g_K\phi_{\star y}(k') \star g_Y(y)$.
\endproof

Protomodularity is stable under the passage to any fiber $Pt_Y\CC$ of the fibration of points, namely to any category whose objects are the split epimorphims above $Y$ an whose morphims are the commutative triangles:	
$$
\xymatrix@=10pt{  
	X \ar[ddrr]_f	\ar[rrrr]^{h} &&&&	X'  \ar[ddll]_{f'} \\
	&&&&&\\
	&&  Y 	\ar@<-1ex>[uull]_<<<<s  \ar@<-1ex>[uurr]_>>>>{s'}
}
$$
Then, in the fiber $Pt_Y\CC$:\\
- the abelian objects are the split epimorphisms $(f,s)$ such that $[R[f],R[f]]=0$; i.e. such that the kernel equivalence relation $R[f]$ (where $xR[f]x'\iff f(x)=f(x')$) is abelian; when "Huq=Smith" is valid (see next section) which is true for $\SkB$ and not for $\DiGp$, they coincide with the split epimorphisms having abelian kernel;\\
- the normal subobjects in $Pt_Y\CC$ are the monomorphism $v$:	
	$$
	\xymatrix@=6pt{  
		\Ker f \ar@{ >->}[dd]_{u} \ar@{ >->}[rr]^{k_f} &&	X \ar@{ >->}[dd]_{v}\ar[rrr]^f	 &&& Y  \ar@<1ex>[lll]^{s} \ar@{=}[dd]\\
		&&&&&\\
		\Ker f' \ar@{ >->}[rr]_{k_{f'}} &&	X' \ar[rrr]_{f'} &&&  Y	 \ar@<-1ex>[lll]_{s'}
	}
	$$
such that $v.k_f=k_{f'}.u$ is normal in $\CC$, see  \cite{strong}.	

Now, there is a result in \cite{Gers} following which:\\
\emph{in the category $\Gp$ of groups, given any diagram of exact sequences:
	$$
	\xymatrix@=6pt{  
		\Ker f \ar@{ >->}[dd]_{u} \ar@{ >->}[rr]^{k_f} &&	X \ar@{ >->}[dd]_{v}\ar@{->>}[rrr]^f	 &&& Y   \ar@{=}[dd]\\
		&&&&&\\
		\Ker f' \ar@{ >->}[rr]_{k_{f'}} &&	X' \ar@{->>}[rrr]_{f'} &&&  Y	 
	}
	$$
	when $u$ is a normal subgroup, so is $k_{f'}.u=v.k_f$.}
	
Accordingly, in the category $\Gp$, for any group $Y$, the (left exact) kernel functor $0^*_Y: Pt_Y\Gp\to \Gp$ reflects normal subobjects. Whence the following natural \cite{strong}:
\begin{definition}
		A pointed category $\CC$ is called strongly protomodular when it is protomodular and such that for any $Y\in\CC$, the kernel functor $0^*_Y: Pt_Y\CC\to \CC$ reflects normal subobjects. 	
	\end{definition} 

Besides $\Gp$, the category $\Rng$ of rings and any category of $R$-algebras are strongly protomodular. However 	the category $\DiGp$ is not strongly protomodular \cite{strong}, it was precisely introduced for that in \cite{strong}.
\begin{theorem}
		The category $\SkB$ is strongly protomodular.
	\end{theorem} 	
\proof
Suppose $u$ normal in $\SkB\;\;$ :
$$
\xymatrix@=6pt{  
	\Ker f \ar@{ >->}[dd]_{u} \ar@{ >->}[rrr]^{k_f} &&&	X \ar@{ >->}[dd]_{v}\ar[rrr]^f	 &&& Y  \ar@<1ex>[lll]^{s} \ar@{=}[dd]\\
	&&&&&\\
	\Ker f' \ar@{ >->}[rrr]_{k_{f'}} &&&	X' \ar[rrr]^{f'} &&&  Y	 \ar@<1ex>[lll]^{s'}
}
$$
then $u:(Kerf,*)\into (Kerf',*)$ and $u:(Kerf,\circ)\into (Kerf',\circ)$ are normal in $\Gp$. So, by Gerstenhaber, both $k_{f'}.u:(Kerf,*)\into (X',*)$ and $k_{f'}.u:(Kerf,\circ)\into (X',\circ)$ are both normal in $\Gp$. It remains to show that $\lambda_{x'}(Kerf)\subset \Ker f$ for any $x'\in X'$.
	
Any $x'\in X'$ is such that $x'=s(y)\circ k'$ for some $k'\in \Ker f'$; so for any $k\in \Ker f$: $\lambda_{x'}(k)=\lambda_{s(y)\circ k'}(k)=\lambda_{s(y)}\lambda_{k'}(k)$. Then:\\
- $\bar k=\lambda_{k'}(k)\in \Ker f$, since $u$ is normal in $SkB$;\\
- $\lambda_{s(y)}(\bar k)\in \Ker f$, since $k_f$ is normal in $SkB$ and $s(y)\in X$.
\endproof

\subsection{"Huq=Smith"}
We recalled that in a pointed protomodular category, there is an intrinsic notion of commutation of subobjects $[u,v]=0$ (Huq commutation) \cite{Huq}. As a Mal'tsev category, there is, as well, an intrinsic notion of commutation of pair of equivalence relations $[R,S]=0$ (Smith commutation) \cite{Smith}. On the other hand, to any equivalence relation $R\stackrel{(d_0,d_1)}{\rightarrowtail}X\times X$ we can associate a normal subobject $I_R=d_1.\ker d_0$ (its normalization):
$$
	\xymatrix@=6pt{  
		\Ker d_0 \ar[dd]	\ar@{ >->}[rrr]_{\ker d_0}  \ar@<2ex>@{ >->}[rrrrrr]^{I_R} &&&	R \ar[dd]_{d_0} \ar[rrr]_{d_1} &&& X\\
		&&&&&\\
		1 \ar@{ >->}[rrr]^{}&&& X 	  
	}
	$$ 	
In any protomodular category,  we have always $[R,S]=0\implies [I_R,I_S]=0$ ("Smith implies Huq"). The converse is not true in general: again, the category $\DiGp$  produces a counterexample. But the converse is true for $\Gp$. It is true for $\SkB$ as well, see \cite{BFP}. In this case, we say that "Huq=Smith" holds. Actually, as previously shown in \cite{BG}, this is true in any strongly protomodular category \cite{BG}. Under "Huq=Smith", in any finitely cocomplete exact pointed protomodular category, the commutator of two equivalence relations coincides with the commutator of their associated normal subobjects, see \cite{strN}, and \cite{BFP} for the specific case of $\SkB$.
	
A quick application of "Huq=Smith" is the following one. In a protomodular category, on any reflexixe graph,
$$
\xymatrix@=6pt{  
	X_1 \ar@<1ex>[rr]^{d_0} \ar@<-1ex>[rr]_{d_1}\ar@<1ex>[rr]^{d_0} && X_0 \ar[ll]|{s_0} 	  
}
$$
there is at most one groupoid structure. As in any Mal'tsev category, this is the case if and only if $[R[d_0],R[d_1]]=0$. So, in $\SkB$, a reflexive graph is underlying a groupoid structure if and only if $[\Ker d_0,\Ker d_1]=0$, as it is the case in $\Gp$ as shown in \cite{stras}.	

\section{Baer sums of exact sequences of skew braces}

We shall now focus our attention on the construction of the Baer sums in $\SkB$. 

\subsection{Brief recall about Baer sums in $\Gp$}

Consider any exact sequence in $\Gp$ with abelian kernel:
$$ \xymatrix@=25pt{
	1 \ar[r] &	{(A,+)\;}  \ar@{>->}[r]^{k_f} & (X,\star) \ar[dr]_{Inn} \ar@{->>}[r]^{f} & (Y,\star) \ar[d]^{\phi} \ar[r] & 1 \\
	& & & \Aut A
}
$$
The inner automorphims of $(X,\star)$ are stable on the normal subgroup $(A,+)$, whence the homomorphism $Inn$. And we get $Inn.k_f=1$ as soon as $(A,+)$ is abelian. So, this produces a factorization $\phi$, defined by $\phi_y(a)=x\star a \star x^{-\star}$ for any $x$ such that $f(x)=y$. In this way, we get a group action $\phi$, and a split epimorphism $(p_Y,\iota_Y): (Y\times A, \ltimes_\phi) \splito (Y,\star)$. Both $\phi$ and this split epimorphism will be call indifferently the \emph{direction} of the given exact sequence (a homomorphism $f: (X,*)\to (Y,*)$ with abelian kernel being thought as an \emph{affine object} in the slice category $\Gp/Y$). This gives rise to the abelian group structure $Opext(Y,A,\phi)$, see \cite{homology}, on the isomorphic classes of exact sequences with same direction $\phi$, which is obtained by completing with the lower exact sequence the following $3\times 3$ diagram:
$$ \xymatrix@=18pt{
	& 1 \ar[d] & 1 \ar[d] && 1 \ar[d] \\
	1 \ar[r] &	{(A,+)\;} \ar@{ >->}[d]_{(-Id_A,Id_A)}  \ar@{>->}[r]^{Id_A} & (A,+) \ar@{ >->}[d]_{(-k_f,k_{f'})} \ar@{->>}[rr]^{} && {1\;} \ar@{ >->}[d]^{} \ar[r] & 1\\
	1 \ar[r] &	{(A\times A,+)\;} \ar@{->>}[d]_{+} \ar@{>->}[r]^{k_f\times k_{f'}} & (X\times_YX',\star) \ar@{.>>}[d]^{\theta} \ar@{->>}[rr]^{f\times_Yf'} && (Y,\star) \ar@{->>}[d]^{Id_Y} \ar[r] & 1\\
	1 \ar[r] &	{(A,+)\;} \ar@{ >.>}[r]_{k_{f\otimes_Y f'}} \ar[d]  & (X\otimes_{_Y}X',\star)   \ar@{.>>}[rr]_{f\otimes_Yf'} \ar[d] && (Y,\star) \ar[r] \ar[d] & 1 \\
	& 1 & 1 && 1 
}
$$
where $f\times_Y f'$ is given  by the following pullback in $\Gp$:
$$
\xymatrix@=6pt{  
	X\times_Y X' \ar@{->>}[ddd]_{p}\ar@{->>}[rrr]^{p'} \ar@{.>>}[rrrddd]^{f\times_Yf'}	 &&& X' \ar@{->>}[ddd]^{f'}\\
	&&&&&\\
	&&&&&\\
	X \ar@{->>}[rrr]_{f} &&&  Y	 
}
$$
namely by the product in the slice category $\Gp/Y$.

\subsection{Baer sums in $\SkB$}

As any pointed protomodular category, the category $\SkB$ has Baer sums of exact sequences with abelian kernel relation and same direction, see \cite{Baer}. Under the "Huq=Smith" condition, the description is greatly simplified.

\begin{proposition}
Let be given any exact sequence with abelian kernel in $\SkB$:
$$
\xymatrix@=7pt{ 
	1 \ar[rr] && (A,+,+) \ar@{ >->}[rr]^k && (X,\star,\circ)  \ar@{->>}[rr]^>>>>{f}  && (Y,\star,\circ) \ar[rr] && 1
}$$
1) $\lambda^X_a(b)=b$ for all $(a,b)\in A\times A$;\\
2) $\lambda^X_x(a)=\lambda^X_{x'}(a)$, for any $(x,x',a)\in R[f]\times A$;\\
3) $u\star v^{-\star}\star w=u\circ v^{-\circ}\circ w$ for any triple $(uR[f]vR[f]w)$;\\
4) given any $(y,a)\in Y\times A$, the mapping $a\mapsto  (a\circ x)\star x^{-\star}$ with $f(x)=y$ is independent from the choice of $x$. We shall denote it by $\xi^{f}_y: A\to A$. 
\end{proposition}
\proof
Condition 1) is straightforward since $(A,+)$ is abelian; actually condition 1) is already true in $\DiGp$. We have $\lambda_x(a)=\lambda_{x'}(a)$ if and only if $\lambda_{x'}(\lambda_x)^{-1}(a)=a$. By the major result $(2)$ in $\SkB$ recalled in the introduction, it is equivalent to $\lambda_{x'\circ x^{-\circ}}(a)=a$, with $x'\circ x^{-\circ}\in A$. Condition 2) is then straightforward from 1).

Now $u\star v^{-\star}\star w=u\circ v^{-\circ}\circ w$ if and only if $v^{-\star}\star w=\lambda^X_u(v^{-\circ}\circ w)$ if and only if $v^{-\circ}\circ w=\lambda^X_{u^{-\circ}}(v^{-\star}\star w)$. This equation is equivalent to:
$\lambda^X_v(v^{-\circ}\circ w)=\lambda^X_v\lambda^X_{u^{-\circ}}(v^{-\star}\star w)=\lambda^X_{v\circ u^{-\circ}}(v^{-\star}\star w)$. Since we have $(uR[f]vR[f]w)$, the second term is $(v^{-\star}\star w)$ by 1), while  the first one is $(v^{-\star}\star w)$ as well, in any case.

Suppose $f(x)=y=f(x')$, then $x'=x\star a$ with $a\in A$; then $((a\circ x^{-\circ})R[f]xR[f]x')$ and: $a\circ x'=(a\circ x)\circ x^{-\circ}\circ x'=(a\circ x)\star x^{-\star}\star x'$. Accordingly: $(a\circ x)\star x^{-\star}=(a\circ x')\star x'^{-\star}$. 
\endproof
The proof of the point 3) precisely shows how, in $\SkB$, any surjection with abelian kernel actually gets an abelian kernel relation (i.e. $[R[f],R[f]]=0$), a conceptual consequence of the "Huq=Smith" condition.
\begin{definition}
	The mapping $\xi^f:Y\to S_A$, defined by the point 4) above, is called the skewing index of the sequence.
\end{definition}
Now define $\lambda^f: Y \to S_A$ by $y\mapsto \lambda^X_x(a)$ with $f(x)=y$ thanks to the point 2) in the above proposition. We then get a non split similar result to Proposition \ref{split}:

\begin{proposition}
	Given any exact sequence of skew brace as above, we get $\xi^f_y\phi_{\circ y}=\phi_{\star y}\lambda^f_y$ for any $y\in Y$.
\end{proposition}
\proof
Choose $x$ such that $f(x)=y$. Then: $\xi^f_y\phi_{\circ y}(a)=((x\circ a \circ x^{-\circ})\circ x)\star x^{-\star}=(x\circ a)\star x^{-\star}$; while: $\phi_{\star y}\lambda^f_y(a)=x\star (x^{-\star}\star(x\circ a))\star x^{-\star}=(x\circ a)\star x^{-\star}$
\endproof

Given any exact sequence in $\SkB$ as above, since 
$$
\xymatrix@=7pt{ 
	1 \ar[rr] && (A,+) \ar@{ >->}[rr]^{k_f} && (X,\star)  \ar@{->>}[rr]^>>>>{f}  && (Y,\star) \ar[rr] && 1\\
		1 \ar[rr] && (A,+) \ar@{ >->}[rr]^{k_f} && (X,\circ)  \ar@{->>}[rr]^>>>>{f}  && (Y,\circ) \ar[rr] && 1
}$$
are exact sequences in $\Gp$ with abelian kernel $(A,+)$, they define, as we recalled in the previous section, two group actions $\phi_\star: (Y,\star) \to \Aut(A,+)$ and $\phi_\circ: (Y,\circ) \to \Aut(A,+)$.

We are now going to define the Baer sum of exact sequences with abelian kernel $(A,+,+)$ in $\SkB$ when they have same pair of group actions  $(\phi_{\star},\phi_{\circ})$ and same skewing index $\xi$. We shall call the triple $(\phi_{\star},\phi_{\circ},\xi)$, the \emph{direction} of the sequence. So, let us consider now a pair of exact sequences in $\SkB$ with same direction:
$$
\xymatrix@=7pt{ 
1 \ar[rr] && (A,+,+) \ar@{ >->}[rr]^{k_f} && (X,\star,\circ)  \ar@{->>}[rr]^>>>>{f}  && (Y,\star,\circ)  \ar[rr] && 1\\
1 \ar[rr] && (A,+,+) \ar@{ >->}[rr]_{k_{f'}} && (X',\star,\circ)  \ar@{->>}[rr]_>>>>{f'}  && (Y,\star,\circ)  \ar[rr] && 1
}$$
and the induced exact sequence on the pullback $X\times_Y X'$ of $f$ and $f'$ in $\SkB$:
$$
\xymatrix@=7pt{ 
	1 \ar[rr] && (A\times A,+,+) \ar@{ >->}[rr]^{k_f\times k_{f'}} && (X\times_YX',\star,\circ)  \ar@{->>}[rr]^>>>>{f\times_Yf'}  && Y \ar[rr] && 1
}$$
\begin{proposition}
	The subobject $(-k_f,k_{f'}): (A,+,+)\into (X\times_YX',\star,\circ)$ defined by $a\mapsto (-a,a)$ is an ideal of the skew brace $(X\times_YX',\star,\circ)$.
\end{proposition}
\proof
Clearly $(A,+)\into (X\times_YX',\star)$ and $(A,+)\into (X\times_YX',\circ)$ are normal in $\Gp$. It remains to check that $\lambda^{X\times_YX'}_{(x,x')}(-a,a)=(-\lambda^X_x(a),\lambda^{X'}_{x'}(a))$ belongs to this subobject, namely to check that  $\lambda^X_x(a)=\lambda^{X'}_{x'}(a)$ when $(x,x')$ belongs to $X\times_YX'$, namely when $f(x)=f'(x')$. This certainly holds since, by the previous proposition,  two exact sequences with same direction has necessarily same $\lambda^f: Y\to S_K$, which means that $\lambda^{X}_{x}(a)=\lambda^{X'}_{x'}(a)$
as soon as $f(x)=y=f'(x')$. 
\endproof

Since the subobject $(-k_f,k_{f'})$ is an ideal, we can complete, on the model of what is done in $\Gp$, the following $3\times 3$ diagram in $\SkB$:
$$ \xymatrix@=18pt{
	& 1 \ar[d] & 1 \ar[d] && 1 \ar[d] \\
 1 \ar[r] &	{(A,+,+)\;} \ar@{ >->}[d]_{(-Id_A,Id_A)}  \ar@{>->}[r]^{Id_A} & (A,+,+) \ar@{ >->}[d]_{(-k_f,k_{f'})} \ar@{->>}[rr]^{} && {1\;} \ar@{ >->}[d]^{} \ar[r] & 1\\
 1 \ar[r] &	{(A\times A,+,+)\;} \ar@{->>}[d]_{+} \ar@{>->}[r]^{k_f\times k_{f'}} & (X\times_YX',\star,\circ) \ar@{.>>}[d]^{\theta} \ar@{->>}[rr]^{f\times_Yf'} && (Y,\star,\circ)  \ar@{->>}[d]^{Id_Y} \ar[r] & 1\\
 1 \ar[r] &	{(A,+,+)\;} \ar@{ >.>}[r]_{k_{f\otimes_Y f'}} \ar[d]  & X\otimes_{_Y}X'   \ar@{.>>}[rr]_{f\otimes_Yf'} \ar[d] && (Y,\star,\circ) \ar[r] \ar[d] & 1 \\
 	& 1 & 1 && 1 
}
$$ 
It produces the lower exact sequence with abelian kernel in $\SkB$. It remains to check it has same direction as the two original ones. It has the same $\star$ and $\circ$ group actions by the Baer constructions in  $\Gp$. It remains to check it has same skewing index:
\proof
The index $\xi_y$ of the lower exact sequence is defined by:\\ $k_{f\otimes_Y f'}(a)\mapsto (k_{f\otimes_Y f'}(a)\circ \theta(x,x'))\star \theta(x,x')^{-\star}$ for any $(x,x')\in X\times_Y X'$ such that $(f\times_Yf')(x,x')=y$, namely $f(x)=y.=f'(x')$. Now:\\ $(k_{f\otimes_Y f'}(a)\circ \theta(x,x'))\star \theta(x,x')^{-\star}=(\theta(0,a)\circ \theta(x,x'))\star \theta(x,x')^{-\star}$\\ =$\theta (0\circ x,a \circ x')\star \theta (x^{-\star},x'^{-\star})$ =$\theta(0,(a \circ x')\star x'^{-\star}))=\theta(0,\xi^{f'}_y(a))$\\=$k_{f\otimes_Y f'}(0+\xi^{f'}_y(a))=k_{f\otimes_Y f'}(\xi^{f'}_y(a))$. This means: $\xi_y(a)=\xi^{f'}_y(a)(=\xi^{f}_y(a))$.
\endproof 

As in $\Gp$, the inverse of the exact sequence:
$$
\xymatrix@=7pt{ 
	1 \ar[rr] && (A,+,+) \ar@{ >->}[rr]^k && (X,\star,\circ)  \ar@{->>}[rr]^>>>>{f}  && (Y,\star,\circ) \ar[rr] && 1
}$$
is  given by the following one:
$$
\xymatrix@=7pt{ 
	1 \ar[rr] && (A,+,+) \ar@{ >->}[rr]^{-k} && (X,\star,\circ)  \ar@{->>}[rr]^>>>>{f}  && (Y,\star,\circ) \ar[rr] && 1
}$$
The unit of this group is given by the sum of the two previous sequences; it is the unique split epimorphism obtained from the direction $(\phi_*,\phi_\circ, \xi)$ by Proposition \ref{unit}. It belongs to $\SkB$ since, by the previous $3\times 3$ lemma construction, its domain is the quotient of $(X\times_YX,*,\circ)$ which belongs to $\SkB$.

According to the notations of Mac Lane \cite{homology}, we shall denote this abelian group by $Opext\bigl((Y,*,\circ),(A,+,+),(\phi_*,\phi_\circ, \xi)\bigr)$. It is clear that:\\ 
$Opext\bigl((Y,*,*),(A,+,+),(\phi_*,\phi_*,\check 1)\bigr)=Opext\bigl((Y,*),(A,+),\phi_*\bigr)$,\\ where $\check 1$ is the trivial index.

\end{document}